\documentclass[12pt]{article}
\newcommand{\be}{\begin{eqnarray}}
\newcommand{\ee}{\end{eqnarray}}
\newcommand{\bee}{\begin{eqnarray*}}
\newcommand{\eee}{\end{eqnarray*}}
\newcommand{\tS}{\widetilde{S}}
\newcommand{\hS}{\widehat{S}}
\newcommand{\hW}{\widehat{W}}
\newcommand{\cF}{\mathcal{F}}
\newtheorem{teo}{Theorem}

\begin{document}

\title{A Delayed Black and Scholes Formula I  }
 % Enter your title between curly braces
\author{By Mercedes Arriojas\thanks {The research of this author is supported
 in part by the Center of Scientific and Human Development  
 in Venezuela, and by Southern Illinois University.}\,,  
Yaozhong Hu\thanks {The research of this author is supported in part by the
 National Science Foundation under Grant No. DMS 0204613 and No. EPS-9874732,
 matching support from the State of Kansas and General Research Fund of the
 University of Kansas.}\,,  
Salah-Eldin Mohammed\thanks {The research of this author is supported in part
 by NSF Grants DMS-9975462 and DMS-0203368, and by the University of Central
 Venezuela}\,,  \\
and Gyula Pap\thanks {The research of this author is supported in part by a 
 Fulbright fellowship.}}
 
\maketitle

\abstract{In this article we develop an explicit formula for pricing European
 options when the underlying stock price follows a non-linear stochastic
 differential delay equation (sdde).
We believe that the proposed  model is sufficiently flexible to fit real
 market data, and is yet simple enough to allow for a closed-form
 representation of the option price.
Furthermore, the model maintains the no-arbitrage property and the completeness of the market.
The derivation of the option-pricing formula is based on an equivalent
 martingale measure. }

\section{Introduction} % Enter section title between curly braces
 
The {\it Black and Scholes Formula} has been one of the most important
 consequences of the study of continuous time models in finance ([Bac], [C], [Me$_1$], [Me$_2$]). However, the
 fitness of the model has been questioned on the basis of the
 assumption of constant volatility ([Sc], [R]), since empirical
 evidence shows that volatility actually depends on time in a way
 that is not predictable. This is sometimes pointed out as the reason for inaccurate 
%cause of lack of accuracy of 
predictions made by the Black and Scholes
 formula. On the other hand, the need for better ways of understanding the behavior of
 many natural processes has motivated the development of dynamic models of
 these processes that take into consideration the influence of past events on
 the current and future states of the system ([I.N], [Ku], [K.N], [Mo${}_1$],
 [Mo${}_2$], [M.T], [E.\O.S]).
This view is specially appropriate in the study of financial variables, since
 predictions about their evolution take strongly into account the knowledge of
 their past ([H.\O], [S.K]).
%In this framework we will derive a formula for pricing options on stocks with
 %hereditary structure (Theorem 4).

In this paper we consider the effect of the past in the determination of the
 fair price of a call option.
In particular, we assume that the stock price satisfies a stochastic functional
 differential equation (sfde).
We consider call options that can be exercised only at the maturity date, viz.
 \emph{European call options}.
We derive an {\it explicit} formula for the valuation of a European call option
 on a given stock (Theorem 4) (cf. ([B.S], [Me$_1$], [H.R]).
Note that based on the preprint version of the present paper, the logarithmic
 utility of an insider has been computed and the stability of the European call
 option has been proved in [St].

Tests of the classical Black and Scholes model against real market data suggest
 the existence of significant levels of randomness in the volatility of the
 stock price, as manifested in the observed phenomenon of {\it frowns} and
 {\it smiles} ([Bat]).
One of the motivations behind our model for the stock price is to account for
 such volatility in a natural manner, while at the same time maintain an
 explicit formula for the option price.
It is hoped that the parameters of the proposed model will allow enough
 flexibility for a better fit than that of the Black and Scholes model  when
 tested against real market data.

%Although options are a very particular class of financial securities, they
% show the characteristic properties of  more general forms of investment such
% as \emph{contingent claims} or \emph{derivatives}.
International markets for contingent claims have experienced remarkable
 growth in the last thirty years.
This makes the study of option pricing of special interest in the present
 context, since this theory may lead to a general theory of pricing contingent
 claims or derivatives with hereditary structure.

\section{Stochastic delay models for the stock price}

In this section we propose a stochastic delay model for the evolution of the
 stock price.
We prove that the proposed model is feasible.
In Section 3, we formulate and solve the option pricing problem for the model.

Consider a stock whose price at time $t$ is given by a stochastic process
 $S(t)$ satisfying the following stochastic delay  differential equation
 (sdde):
 \begin{equation}
\left.
\begin{array}{lll}
dS(t)& = & \mu S(t-a)S(t)\,dt + g(S(t-b))S(t)\,dW(t),\quad t\in [0,T]
\nonumber \\
S(t)&=&\varphi (t), \quad t\in [-L,0].
\end{array}
\right\}
\end{equation}
on a probability space $(\Omega, {\cal F}, P)$ with a filtration
 $({\cal F}_t)_{0 \leq t \leq T}$ satisfying the usual conditions.
In the above sdde, $\mu, a, b$ are
 positive constants,  $L:=\max \{a,b\}$, and $g:\textbf{R}\to \textbf{R}$ is a continuous
function. The initial process  $\varphi: \Omega \to C([-L,0],\textbf{R})$ is
 ${\cal F}_0$-measurable with respect to the Borel $\sigma$-algebra of
 $C([-L,0],\textbf{R})$.
The process $W$ is a one-dimensional standard Brownian motion adapted to the
 filtration $({\cal F}_t)_{0 \leq t \leq T}$.

In our next result, we will show that the above model  is feasible in the 
sense that it admits a pathwise unique solution $S$ such that
 $S(t)> 0$ almost surely for all $t \geq 0$ whenever 
 $\varphi (0) > 0$.

\medskip
%\medskip

\begin{teo}

The sdde (1) has a pathwise unique solution $S$
 for a given ${\cal F}_0$-measurable initial process
 $\varphi: \Omega \to C([-L,0],\textbf{R})$.
Furthermore, if $\varphi (0) \geq 0$ ($\varphi (0) > 0$) a.s., then
 $S(t) \geq 0$ ($S(t) > 0$) for all $t \geq 0$ a.s..
%If in addition $\varphi (0) > 0$ a.s., then  $S(t) > 0$  for all
% $t \geq 0$ a.s..

\end{teo}

\noindent
\textsl{Proof.}

%\begin{\itemize}

%\item{}
\noindent
%{\it Case 1:} $f=f_1$.

Suppose $\varphi (0) \geq 0$. Define $l:=\min\{ a,b \} > 0$ and let $t\in [0,l]$.  Then (1) gives
 \begin{equation}
\left.
\begin{array}{lll}
dS(t)& = & S(t)[\mu \varphi(t-a)\,dt + g(\varphi (t-b))\,dW(t)],
\quad t\in [0,l] \nonumber \\
S(0)&=&\varphi (0).
\end{array}
\right\}
\end{equation}
Define the semimartingale
$$ N(t):= \mu \int_0^t \varphi(u-a)\,du + \int_0^t g(\varphi (u-b))\,dW(u),
\quad t \in [0,l],$$
 and denote by $[N,N](t)= \int_0^t g(\varphi(u-b))^2\,du, \, t \in [0,l]$, its
 quadratic variation.
Then (2) becomes
%
%\begin{equation}
$$dS(t)= S(t)\, dN(t), \quad t > 0, \qquad   S(0)=\varphi (0), $$
% \end{equation}
%
which has the unique solution
%of (5) is given by
%
 \begin{eqnarray*}
S(t)&=& \varphi (0) \exp \{N(t)-\frac{1}{2}[N,N](t)\},\\
&=& \varphi (0) \exp\Big \{\mu \int_0^t \varphi(u-a)\,du\\
&&\quad 
       + \int_0^t g(\varphi (u-b))\,dW(u)
       -\frac{1}{2} \int_0^t g(\varphi(u-b))^2\,du\Big \},
\end{eqnarray*}
for $t \in [0,l]$.
This clearly implies that $S(t)> 0$ for all $t \in [0,l]$ almost surely,
when $\varphi (0) > 0$ a.s.
By a similar argument, it follows that $S(t) > 0$  for all $t \in [l,2l]$ a.s..
Therefore $S(t) > 0$  for all $t \geq 0$ a.s., by induction.
Note that the above argument also gives existence and pathwise-uniqueness
 of the solution to (1).

\medskip

\section{A delayed option pricing formula}

% Enter subsection title between curly braces

Consider a market consisting of a riskless asset (e.g., a bond or bank account)
 $B(t)$ with rate of return $r\geq0$ (i.e., $B(t)=e^{rt}$\,) and a single stock
 whose price $S(t)$ at time $t$ satisfies the sdde (1) where $\varphi(0)>0$
 a.s.. 
In the sdde (2), assume further that the delays $a,b$ are positive and $g$ is
 continuous.
Consider an option, written on the stock, with maturity at some future time
 $T > t$ and an exercise price $K$.
Assume also  that there are no transaction costs and that the underlying stock
 pays no dividends.
Our main objective is to derive the fair price of the option at time $t$.
In the following discussion, we will  obtain an equivalent martingale
 measure with the help of Girsanov's theorem.

Let
 $$\tS(t):=\frac{S(t)}{B(t)}=e^{-rt}S(t),\qquad t\in[0,T],$$
 be the discounted stock price process.
Then by It\^o's formula (the product rule), we obtain
 \begin{eqnarray*}
  d\tS(t)&=&e^{-rt}dS(t)+S(t)(-re^{-rt})\,dt\\
         &=&\tS(t)\Big[\big \{\mu S(t-a)-r\big \}\,dt+g\big(S(t-b)\big 
)\,dW(t)\Big].
 \end{eqnarray*}
Let
 $$\hS(t):=\int_0^t \big \{\mu S(u-a)-r\big \}\,du
           +\int_0^tg\big(S(u-b)\big)\,dW(u),\qquad t\in[0,T].
$$
Then
 \begin{equation}\label{tShS}
  d\tS(t)=\tS(t)\,d\hS(t), \quad 0 < t < T.
 \end{equation}
Taking into account that $\tS(0)=\varphi(0)$, we have
 \begin{equation}\label{tS}
  \tS(t)=\varphi(0)+\int_0^t\tS(u)\,d\hS(u),\qquad t\in[0,T].
 \end{equation}

We now recall Girsanov's theorem (see, e.g., Theorem 5.5 in [K.K]).

\begin{teo}

\textbf{\textup{(Girsanov)}}
Let $W(t)$, $t\in[0,T]$, be a standard Wiener process on
 $(\Omega,\mathcal{F},P)$.
Let $\Sigma $ be a predictable process such that
 $\int_0^T|\Sigma (u)|^2du<\infty$ a.s., and let
 $$\varrho_t
   :=\exp\left\{\int_0^t \Sigma (u)\,dW(u)-\frac{1}{2}\int_0^t|\Sigma 
(u)|^2\,du\right\},
   \qquad t\in[0,T].$$
Suppose that $E_P(\varrho_T)=1$, where $E_P$ denotes expectation with respect
 to the probability measure $P$.
Define the probability measure $Q$ on $(\Omega,\mathcal{F})$ by
 $dQ :=\varrho_T\,dP$.
Then the process
 $$\hW(t):=W(t)-\int_0^t \Sigma (u)\,du,\qquad t\in[0,T],$$
 is a standard Wiener process under the measure $Q$.

\end{teo}

\medskip

From now on, we will assume that the function $g:\textbf{R}\to \textbf{R}$ in
 the sdde (1) satisfies the following hypothesis:
\bigskip

\noindent
{\bf{Hypothesis (B).}} \ $g(v)\neq 0$ whenever $v \neq 0$.
%$g:\textbf{R}\to \textbf{R}$ as considered before is such that

\bigskip

We want to apply Girsanov's theorem with the process
 $$\Sigma (u):=-\frac{\{\mu S(u-a)-r\}}{g\big(S(u-b)\big)},\qquad u\in[0,T].$$
Hypothesis (B) implies that $\Sigma $ is well-defined, since by Theorem 1, 
$S(t)>0$
  for all $t\in[0,T]$ a.s..
Clearly $\Sigma (t)$, $t\in[0,T],$ is a predictable process.
Moreover, $\int_0^T|\Sigma (u)|^2du<\infty$ a.s., since sample-path continuity
 of the process $S(t)$, $t\in[0,T]$, implies almost sure boundedness of $S(t)$,
 $t\in[0,T]$, and Hypothesis (B) implies that $1/g(v)$, $v\in(0,\infty)$, is
 bounded on bounded intervals.
Now let $l:=\min(a,b)$. Set \ $\cF_t:=\cF_0$ for $t \leq 0$.
Then $\Sigma (u)$, $u\in[0,T]$, is measurable with respect to the
 $\sigma$-algebra $\cF_{T-l}$.  
Hence, the stochastic integral $\int_{T-l}^T \Sigma (u)\,dW(u)$
 {\it conditioned on} $\mathcal{F}_{T-l}$ has a normal distribution with mean
 zero and variance $\int_{T-l}^T \Sigma (u)^2\, du$.
Consequently, by the formula for the moment generating function of a normal
 distribution, we obtain
 $$
E_P\left(\exp\left\{\int_{T-l}^T \Sigma (u)\,dW(u)\right\}\,\bigg|\,
            \mathcal{F}_{T-l}\right)
   =\exp\left\{\frac{1}{2}\int_{T-l}^T|\Sigma (u)|^2\,du\right\}
   $$
a.s..
Hence
 $$
E_P\left(\exp\left\{\int_{T-l}^T \Sigma (u)\,dW(u)
                     -\frac{1}{2}\int_{T-l}^T|\Sigma (u)|^2du\right\}\,\bigg|\,
            \mathcal{F}_{T-l}\right) =1
 $$
a.s.. Now the above relation easily implies that
\begin{eqnarray*}
E_P\biggl (\exp\biggl \{\int_0^T \Sigma (u)\,dW(u)
                        -\frac{1}{2}\int_0^T|\Sigma (u)|^2du\biggr\}\,
\bigg|\, \mathcal{F}_{T-l}\biggr )\\
 = \exp\biggl\{\int_0^{T-l} \Sigma (u)\,dW(u)
               -\frac{1}{2}\int_0^{T-l} |\Sigma (u)|^2\, du\biggr\}
\end{eqnarray*}
 %\begin{eqnarray*}
%$$E_P\left(\exp\left\{\int_0^T \Sigma (u)\,dW(u)
%                      -\frac{1}{2}\int_0^T|\Sigma (u)|^2du\right\}\,
%\bigg|\, \mathcal{F}_{T-l}\right)$$
% \\
%$$= \exp\left\{\int_0^{T-l} \Sigma (u)\,dW(u)
%               -\frac{1}{2}\int_0^{T-l} |\Sigma (u)|^2\, du\right\}
%$$
%\end{eqnarray*}
%
a.s..  
Let $k$ to be a positive integer such that $0 \leq T-kl \leq l$. 
Then by successive conditioning using  backward steps of length $l$, an
 inductive argument gives
\begin{eqnarray*}
E_P\biggl (\exp\biggl \{\int_0^T \Sigma (u)\,dW(u)
           -  \frac{1}{2}\int_0^T|\Sigma (u)|^2du\biggr \}\,
           \bigg|\, \mathcal{F}_{T-kl}\biggr ) \\
= \exp\biggl \{\int_0^{T-kl} \Sigma (u)\,dW(u)
               -\frac{1}{2}\int_0^{T-kl} |\Sigma (u)|^2\, du\biggr \}
\end{eqnarray*}
a.s.. 
Taking conditional expectation with respect to $\mathcal{F}_0$ on both sides of
 the above equation, we obtain 
\begin{eqnarray*}
E_P\biggl (\exp\biggl \{\int_0^T \Sigma (u)\,dW(u)
           -  \frac{1}{2}\int_0^T|\Sigma (u)|^2du\biggr\}\,
           \bigg|\, \mathcal{F}_0\biggr) \\
=E_P\left(\exp\left\{\int_0^{T-kl} \Sigma (u)\,dW(u)
          -  \frac{1}{2}\int_0^{T-kl}|\Sigma (u)|^2du\right\}\,
          \bigg|\, \mathcal{F}_{0}\right)=1
\end{eqnarray*}
a.s.. 
Taking the expectation of the above equation, we immediately obtain 
 $$
  E_P(\varrho_T)=1.
 $$
Therefore, the Girsanov's theorem (Theorem 2) applies and  the process
 $$
 \hW(t)
  :=W(t)+\int_0^t\frac{\{\mu S(u-a)-r\}}{g\big(S(u-b)\big)}du,\qquad t\in[0,T],
 $$
 is a standard Wiener process under the measure $Q$ defined by
 $dQ :=\varrho_T\,dP$ with
 $$\varrho_T
   :=\exp\left\{-\int_0^T\frac{\{\mu S(u-a)-r\}}{g\big(S(u-b)\big)}\,dW(u)
                -\frac{1}{2}
                 \int_0^T\left|\frac{\mu S(u-a)-r}{g\big(S(u-b)\big)}\right|^2
                  du\right\}
 $$
a.s..
 %for $t\in[0,T]$.
Since the process $\hS(t)$, $t\in[0,T]$, can be written in the form
 \begin{equation}\label{hShW}
  \hS(t)=\int_0^tg\big(S(u-b)\big)\,d\hW(u),\qquad t\in[0,T],
 \end{equation}
 we conclude that $\hS(t)$, $t\in[0,T]$, is a continuous $Q$-martingale (i.e.,
 a continuous martingale under the measure $Q$\,).
Furthermore, by the representation (\ref{tS}), the discounted stock price
 process $\tS(t)$, $t\in[0,T]$, is also a continuous $Q$-martingale.
In other words, $Q$ is an equivalent martingale measure.
By the well-known theorem on trading strategies (e.g., Theorem 7.1 in
 [K.K]), it follows that the market consisting of $\{B(t),S(t):t\in[0,T]\}$
 satisfies the no-arbitrage property: There is no admissible
 self-financing strategy which gives an arbitrage opportunity.

We now establish the completeness of the market 
$\{ B(t), S(t): t \in [0,T] \}$.

From the proof of Theorem 1, it follows that the solution of the sdde
 (1) satisfies the relation
 $$S(t)=\varphi(0)
        \exp\left\{\int_0^tg\big(S(u-b)\big)\,dW(u)
                   +\mu\int_0^tS(u-a)\,du
                   -\frac{1}{2}\int_0^tg\big(S(u-b)\big)^2du\right\}
 $$
 a.s. for $t\in[0,T]$. Hence we have
 \begin{equation}\label{tShW}
  \tS(t)=\varphi(0)
         \exp\left\{\int_0^tg\big(S(u-b)\big)\,d\hW(u)
                    -\frac{1}{2}\int_0^tg\big(S(u-b)\big)^2du\right\}
 \end{equation}
 for $t\in[0,T]$. 
 %[$S(t)$ and $\tilde S(t)$ are the same 
%except we consider $\tilde S(t)$ as a functional of $\tilde W$ and 
%$S(t)$ is considered as a functional of $W$].
By examining the definitions of $\tilde S, \hW, \hS$ and equation (6), it is not hard to see that for $t \geq 0$,
$\cF_t^S=\cF_t^{\tS}=\cF_t^{\hW}=\cF_t^{W}$, 
%This implies that \ $\cF_t^S=\cF_t^{\tS}=\cF_t^{\hW}=$\cF_t^{W}$, 
the $\sigma$-algebras generated by $\{S (u): u \leq t
\}$, $\{\tS (u): u \leq t \}$,
 $\{\hW(u): u \leq t \}$, $\{W(u): u \leq t \}$, respectively. (Clearly,
  $\cF_t^{W} \subseteq \cF_t$.)
Now, let $X$ be a contingent claim, viz. an integrable non-negative $\cF_T^S$-measurable random
 variable.  Consider the $Q$-martingale
 $$M(t):=E_Q(e^{-rT}X\,|\,\cF_t^S)
        =E_Q(e^{-rT}X\,|\,\cF_t^{\hW}),\qquad t\in[0,T].$$
By the martingale representation theorem (e.g., Theorem 9.4 in [K.K]),
 there exists an $(\cF_t^{\hW})$-predictable process $h_0(t)$, $t\in[0,T]$, such
 that
 $$\int_0^T h_0(u)^2\,du<\infty\qquad a.s.,$$
 and
 $$M(t)=E_Q(e^{-rT}X)+\int_0^th_0(u)\,d\hW(u),\qquad t\in[0,T].$$
By (\ref{tShS}) and (\ref{hShW}) we obtain
 $d\tS(t)=\tS(t)g\big(S(t-b)\big)\,d\hW(u), \, \,t\in[0,T] $. Define  
 $$\pi_S(t):=\frac{h_0(t)}{\tS(t)g\big(S(t-b)\big)},\qquad
   \pi_B(t):=M(t)-\pi_S(t)\tS(t),
   \qquad t\in[0,T].$$
Consider the strategy $\{(\pi_B(t),\pi_S(t)):t\in [0,T]\}$ which consists of
 holding $\pi_S(t)$ units of the stock and $\pi_B(t)$ units of the bond at time
 $t$.
The value of the portfolio at any time $t \in [0,T]$ is given by 
 $$V(t):=\pi_B(t)e^{rt}+\pi_S(t)S(t)=e^{rt}M(t).$$
 Therefore, by the product rule and the definition of the strategy 
 $\{(\pi_B(t),\pi_S(t)):t\in [0,T]\}$, it follows that 
 $$
 dV(t)=e^{rt}dM(t)+M(t)de^{rt}
        =\pi_B(t)de^{rt}+\pi_S(t)dS(t), \qquad t\in[0,T].
$$
Consequently, $\{(\pi_B(t),\pi_S(t)):t\in[0,T]\}$ is a self-financing strategy.
Moreover, $V(T)=e^{rT}M(T)=X$ a.s.. Hence the contingent claim $X$ is attainable.
This shows that the market $\{B(t), S(t): t \in [0,T] \}$
%(consisting of the riskless asset and the stock)
 is complete, since every contingent claim is attainable.
Moreover, in order for the augmented market  $\{B(t), S(t), X: t \in [0,T] \}$ 
%(consisting of the riskless asset, the stock and the contingent claim $X$\,) 
 to satisfy the no-arbitrage property, the price 
 %$V(t)$ 
 of the claim $X$ must
 be
 $$V(t)=e^{-r(T-t)}E_Q(X\,|\,\cF_t^S)$$
at each $t \in [0,T]$ a.s.. See, e.g., [B.R] or Theorem 9.2 in [K.K].

\bigskip

The above discussion may be summarized in the following
 formula for the fair price $V(t)$ of an option on the stock whose evolution is
 described by the sdde (1).

\medskip

\begin{teo}

Suppose that the stock price $S$ is given by the sdde (1), where
 $\varphi(0)>0$ and $g$ satisfies Hypothesis (B).
Let $T$ be the maturity time of an option (contingent claim) on the stock with
 payoff function $X$, i.e., $X$ is an $\cF_T^S$-measurable non-negative integrable
 random variable.
Then at any time $t \in [0,T]$, the fair price $V(t)$ of the option is given
 by the formula
\begin{eqnarray}
V(t)=e^{-r(T-t)}E_Q(X \,|\,\cF_t^S),
\end{eqnarray}
where $Q$ denotes the probability measure on $(\Omega,\mathcal{F})$ defined by
 $dQ :=\varrho_T\,dP$ with
 $$
\varrho_t
   :=\exp\left\{-\int_0^t\frac{\{\mu S(u-a)-r\}}{g\big(S(u-b)\big)}\,dW(u)
                -\frac{1}{2}
                 \int_0^t\left|\frac{\mu S(u-a)-r}{g\big(S(u-b)\big)}\right|^2
                  du\right\}
$$
 for $t\in[0,T]$.
The measure $Q$ is a martingale measure and the market is complete.

Moreover, there is an adapted and square integrable process $h_0(u), 
\ u\in[0, T]$ such that
\[
E_Q(e^{-rT}X\,|\,\cF_t^S)=E_Q(e^{-rT}X)+\int_0^th_0(u)\,d\hW(u),\qquad t\in[0,T],
\]
where $\hW$ is a standard $Q$-Wiener process given by 
 $$
 \hW(t)
  :=W(t)+\int_0^t\frac{\{\mu S(u-a)-r\}}{g\big(S(u-b)\big)}du,\qquad t\in[0,T],
 $$
The hedging strategy is given by 
\begin{equation}
 \pi_S(t):=\frac{h_0(t)}{\tS(t)g\big(S(t-b)\big)},\qquad
 \pi_B(t):=M(t)-\pi_S(t)\tS(t),
 \qquad t\in[0,T].
\end{equation}

\end{teo}

\medskip
 The following result is a consequence of Theorem 3. It gives a
Black-Scholes-type formula for the value of 
a European option on the stock at any time prior to maturity.

\medskip

\begin{teo}

Assume the conditions of Theorem 3. Let $V(t)$ be the fair price of a European 
 call option written on the stock $S$ with exercise price $K$ and maturity
 time $T$. 
Let $\Phi$ denote the distribution function of the standard normal law,
 i.e.,
 $$
  \Phi(x)
  :=\frac{1}{\sqrt{2\pi}} \int_{-\infty}^x e^{-u^2/2}\, du,
  \qquad x \in \textbf{R}.$$
 Then for all $t \in [T-l,T]$ (where $l:=\min\{a,b\}$), $V(t)$ is given by 
%
 %$$
\begin{eqnarray}\label{VBS}
V(t)=S(t)\Phi(\beta_+(t))-Ke^{-r(T-t)}\Phi(\beta_-(t)), 
%\quad t \in [T-l,T]
\end{eqnarray}
%$$
%
 where
 $$\beta_\pm(t):=\frac{\log\frac{S(t)}{K}
                    +\int_t^T\left(r\pm\frac{1}{2}g(S(u-b))^2\right)du}
                   {\sqrt{\int_t^Tg(S(u-b))^2du}}.
$$
If $T > l$ and $t<T-l$, then 
 \begin{eqnarray}\label{VBSC}
V(t)  & =& e^{rt}
    E_Q\left(H\left(\tS(T-l),-\frac{1}{2}\int_{T-l}^Tg\big(S(u-b)\big)^2du,
\right.\right.
\nonumber\\
&&\left. \left. 
                     \int_{T-l}^Tg\big(S(u-b)\big)^2du\right)
             \,\bigg|\,\cF_t\right)
\end{eqnarray}
 where $H$ is given by 
 $$
H(x,m,\sigma^2):=x e^{m + \sigma^2/2} \Phi(\alpha_1(x,m,\sigma))
                   -Ke^{-rT}\Phi(\alpha_2(x,m,\sigma)),
$$
 and
 $$\alpha_1 (x,m,\sigma):=\frac{1}{\sigma} \biggl [\log
\left(\frac{x}{K}\right) +rT +m+ \sigma^2 
\biggr ],$$
$$\alpha_2 (x,m,\sigma):=\frac{1}{\sigma} \biggl [\log
\left(\frac{x}{K}\right) +rT +m \biggr 
],$$
for $\sigma, x \in \textbf{R}^+, \, m \in \textbf{R}$.

The hedging strategy is given by
\[
 \pi_S(t)=\Phi(\beta_+(t)),\qquad
 \pi_B(t)=-Ke^{-rT}\Phi(\beta_-(t)),\qquad
 t\in[T-\ell,T].
\]
\end{teo}
\bigskip
\noindent
{\bf{Remark 2.}}

 \medskip
 %\begin{itemize}
%\item[(i)] 
If $g(x)=1$ for all $x\in\textbf{R}^+$ then equation (\ref{VBS}) reduces to the
 classical Black
 and Scholes formula. Note that, in contrast with the classical (non-delayed) 
Black and Scholes formula, the fair
 price $V(t)$  in a general delayed model considered in Theorem 4
 depends not only on the stock
 price $S(t)$ at the present time $t$, but also on the whole segment 
$\{S(v):v\in[t-b,T-b]\}$.
(Of course $[t-b,T-b]\subset[0,t]$ since $t\geq T-l$ and $l\leq b$.)
%
%Consequently, one can not obtain the above formula for the fair price $V(t)$
% by the original method based on partial differential equations (as presented,
 %e.g., in Chapter 10 in [K.K]).

%\pagebreak
 
%\end{itemize}
\medskip
\bigskip
\noindent
\textsl{Proof of Theorem 4.}

\medskip

Consider a European call option in the above market with exercise price $K$ and
 maturity time $T$. 
Taking $X= (S(T)-K)^+$ in Theorem 3, the fair price $V(t)$ of the option is
 given by
 $$
  V(t)=e^{-r(T-t)}E_Q(\big(S(T)-K\big)^+\,|\,\cF_t)
      =e^{rt}E_Q(\big(\tS(T)-Ke^{-rT}\big)^+\,|\,\cF_t),
 $$
 at any time $t \in [0,T]$.

We now derive an explicit formula for the option price $V(t)$ at any time  
 $t\in[T-l,T]$.
The representation (\ref{tShW}) of $\tS(t)$ implies that  
 $$
  \tS(T)=\tS(t)
         \exp\left\{\int_t^Tg\big(S(u-b)\big)\,d\hW(u)
                    -\frac{1}{2}\int_t^Tg\big(S(u-b)\big)^2du\right\}
 $$
 for all $t\in[0,T]$. 
Clearly $\tS(t)$ is $\cF_t$-measurable.
If $t\in[T-l,T]$, then $-\frac{1}{2}\int_t^Tg\big(S(u-b)\big)^2du$ is also
 $\cF_t$-measurable.   
Furthermore, when conditioned on $\cF_t$, the distribution of
 $\int_t^Tg\big(S(u-b)\big)\,d\hW(u)$ under $Q$ is the same as that of
 $\sigma\xi$, where $\xi$ is a Gaussian $N(0,1)$-distributed random variable,
 and $\sigma^2=\int_t^Tg\big(S(u-b)\big)^2du$.
Consequently, the fair price at time $t$ is given by
 $$
 V(t)=e^{rt}H\left(\tS(t),-\frac{1}{2}\int_t^Tg\big(S(u-b)\big)^2du,
                     \int_t^Tg\big(S(u-b)\big)^2du\right),
 $$
 where
 $$H(x,m,\sigma^2):=E_Q (xe^{m+\sigma\xi}-Ke^{-rT})^+, \quad \sigma, x \in 
\textbf{R}^+, \, m \in \textbf{R}.$$
Now, an elementary computation yields the following:
 $$
H(x,m,\sigma^2)=x e^{m + \sigma^2/2} \Phi(\alpha_1(x,m,\sigma))
                   -Ke^{-rT}\Phi(\alpha_2(x,m,\sigma)).
$$

Therefore, $V(t)$ takes the form:
 $$V(t)=S(t)\Phi(\beta_+)-Ke^{-r(T-t)}\Phi(\beta_-),$$
 where
 $$\beta_\pm=\frac{\log\frac{S(t)}{K}
                    +\int_t^T\left(r\pm\frac{1}{2}g(S(u-b))^2\right)du}
                   {\sqrt{\int_t^Tg(S(u-b))^2du}}.
$$
For $T>l$ and $t<T-l$, from the representation (\ref{tShW}) of $\tS(t)$, we
 have
 $$\tS(T)=\tS(T-l)
         \exp\left\{\int_{T-l}^Tg\big(S(u-b)\big)\,d\hW(u)
                    -\frac{1}{2}\int_{T-l}^Tg\big(S(u-b)\big)^2du\right\}.$$
Consequently, the option price at time $t$ with $t<T-l$ is given by
 $$
V(t)   =e^{rt}
    E_Q\left(H\left(\tS(T-l),-\frac{1}{2}\int_{T-l}^Tg\big(S(u-b)\big)^2du,
                     \int_{T-l}^Tg\big(S(u-b)\big)^2du\right)
             \,\bigg|\,\cF_t\right).
$$

To calculate the hedging strategy for $t\in[T-\ell,T]$, it suffices to use an
 idea from [B.R], pages 95--96.  
This  completes the proof of the theorem.  $\diamond$

%\bigskip
 
%\newpage
\bigskip
\noindent
{\bf{Remark 4.}}
%\textsl{Remark.}
 \medskip

In the last delay period $[T-l,T]$, one can rewrite the option price
 $V(t), t \in [T-l,T]$ in terms of the solution of a random Black-Scholes pde
 of the form
 \begin{equation}
\left.
\begin{array}{lll}
\displaystyle \frac{\partial F(t,x)}{\partial t}&=& -\displaystyle\frac{1}{2} 
g(S(t-b))^2 x^2 \frac{ \partial^2 F(t,x)}{\partial x^2} -
rx\frac{ \partial F(t,x)}{\partial x} +
rF(t,x), \quad 0 < t <T
\nonumber \\
F(T,x)&=&(x-K)^+, \quad x > 0.
\end{array} \right\}
\end{equation}
The above time-dependent random final-value problem admits a unique
 $(\cF_t)_{t \geq 0}$-adapted random field $F(t,x)$. 
Using the classical It\^o-Ventzell formula ([Kun]) and (7) of Theorem 3, it
 can be shown that 
$$ V(t)= e^{-r(T-t)} F(t, S(t)), \quad t \in [T-b,T].$$
Note that the above representation is no longer valid if $ t \leq T-b$, because
 in this range, the solution $F$ of the final-value problem (11) is
 {\it anticipating} with respect to the filtration $(\cF_t)_{t \geq 0}$.

%\bigskip

%\bigskip
\bigskip

\noindent
{\bf{Acknowledgments.}}

\medskip

The authors are very grateful to R. Kuske and B. $\emptyset$ksendal  
%and the referees 
for very useful suggestions and
corrections to earlier versions of the manuscript. The authors also
acknowledge helpful comments and discussions with Saul Jacka.

%\bigskip
%\bigskip
\bigskip
%\newpage

\centerline{\bf{ References}}

\medskip

\medskip
\par\noindent
[Bac] Bachelier, L.,
Theorie de la Speculation,
{\it Annales de l'Ecole Normale Superieure} 3,
Gauthier-Villards, Paris (1900).

\medskip
\par\noindent
[Bat] Bates, D. S.,
Testing option pricing models,
{\it Statistical Models in Finance}, Handbook of Statistics 14, 567-611,
North-Holland, Amsterdam (1996).

\medskip
\par\noindent
[B.R] Baxter, M. and Rennie, A.,
{\it Financial Calculus}, Cambridge University Press (1996).

\medskip
\par\noindent
[B.S] Black, F. and  Scholes, M.,
The pricing of options and corporate liabilities,
{\it Journal of Political Economy} 81 (May-June 1973), 637-654.

\medskip
\par\noindent
[C] Cootner, P. H.,
{\it The Random Character of Stock Market Prices},
MIT Press, Cambridge, MA (1964).

%\medskip
%\par\noindent
%[D] Davis, M. H. A.,
%Two quick derivations of the Black-Scholes option pricing formula,
%Preprint 91-13, Mathematics Institute, University of Oslo, 1991.

\medskip
\par\noindent
[E.\O.S] Elsanousi, I., \O ksendal, B., and Sulem, A.,
Some solvable stochastic control problems with delay,
{\it Stochastics Stochastics Rep.} 71, no. 1-2, (2000), 69--89.
%{\it Stochastics and Stochastic Reports} 71, no. 1-2, (2000), 69-89. 

%\medskip
%\par\noindent
%[H] Hida, T.,{\it Brownian Motion},
%Springer -Verlag, New York-Heidelberg-Berlin (1980).

\medskip
\par\noindent
[H.R]  Hobson, D.,  and Rogers, L. C. G.,
Complete markets with stochastic volatility,
{\it  Math. Finance\/} 8 (1998), 27--48.

%\medskip
%\par\noindent
%[H] Hu, Y. It\^o-Wiener chaos expansion with exact residual and
% correlation, variance inequalities, {\it J. Theoret. Probab.}, 10 (1997), 
% 835--848.

%\medskip
%\par\noindent
%[H.M.Y] Hu, Y.,  Mohammed, S.-E. A., and Yan, F., 
%Discrete-time approximations of  stochastic delay equations: The 
%Milstein scheme, {\it The Annals of  Probability}, Vol. 32, No. 1A, (2004), 265-314.

\medskip
\par\noindent
[H.\O] Hu, Y. and \O ksendal, B.,  Fractional white noise calculus and
applications to finance, {\it Infinite Dimensional Analysis, Quantum Probability
and Related Topics}, 6 (2003), 1-32.

\medskip
\par\noindent
[I.N]  It\^o, K.,  and Nisio, M.,
On stationary solutions of a stochastic differential equation,
{\it J. Math. Kyoto University\/} 41 (1964), 1--75.

\medskip
\par\noindent
[K.K] Kallianpur, G. and Karandikar R. J.,
{\it Introduction to Option Pricing Theory},
Birkh\"{a}user, Boston-Basel-Berlin (2000).

%\medskip
%\par\noindent
%[K.P] Kloeden, P. E. and Platen E.,
%{\it Numerical Solutions of Stochastic Differential Equations},
%Springer-Verlag, New York-Berlin (1995).

%\medskip
%\par\noindent
%[K.P.S] Kloeden, P. E., Platen E., and Schurz, H.,
%{\it Numerical Solutions of SDE Through Computer Experiment},
%Springer-Verlag, New York-Berlin (1997).

\medskip
\par\noindent
[K.N] Kolmanovskii, V. B., and Nosov, V. R.,
{\it Stability of Functional Differential Equations,},
Academic Press, London, Orlando (1986).

%\medskip
%\par\noindent
%[K.S] Karatzas, I. and Shreve S.,
%{\it Brownian Motion and Stochastic Calculus},
%Springer-Verlag, New York-Berlin (1987).

\medskip
\par\noindent
[Ku] Kushner H. J.,
On the stability of processes defined by stochastic differential-difference
 equations, {\it J. Diff. Equations,} 4 (1968), 424-443.

%\newpage
\medskip
\par\noindent
[Kun]  Kunita, H., {\it Stochastic Flows and Stochastic Differential Equations,}
Cambridge University  Press, Cambridge, New York, Melbourne, Sydney (1990).

%\medskip
%\par\noindent
%[Ma] Malliaris, A. G.,
%It\^{o}'s calculus in financial decision making,
%{\it SIAM Review}, Vol 25, 4, October (1983).

%\medskip
%\par\noindent
%[Man] Mania, M.,
%Derivation of a Generalized Black-Scholes Formula,
%{\it Proceedings of A. Razmadze Mathematical Institute},
%Vol. 115 (1997), 121-148.

\medskip
\par\noindent
[Me${}_1$] Merton, R. C.,
Theory of rational option pricing,
{\it Bell Journal of Economics and Management Science} 4,
Spring (1973), 141-183.

\medskip
\par\noindent
[Me${}_2$] Merton, R. C.,
{\it Continuous-Time Finance}, Basil Blackwell, Oxford (1990).

\medskip
\par\noindent
[M.T] Mizel, V. J. and Trutzer, V.,
Stochastic hereditary equations: existence and asymptotic stability,
{\it Journal of Integral Equations\/} 7 (1984), 1--72.

\medskip
\par\noindent
[Mo${}_1$] Mohammed, S.-E. A.,
{\it Stochastic Functional Differential Equations.} Pitman 99 (1984).

\medskip
\par\noindent
[Mo${}_2$] Mohammed, S.-E. A.,
Stochastic differential systems with memory:
Theory, examples and applications. In ``Stochastic Analysis",
Decreusefond L. Gjerde J., \O ksendal B., Ustunel A.S.,  edit.,
{\it Progress in Probability} 42, Birkhauser (1998), 1-77.

%\medskip
%\par\noindent
%[\O] \O ksendal, B., {\it Stochastic Differential Equations}, Springer, 
%fifth edition (1998).

\medskip
\par\noindent
[R] Rubinstein, M., Implied binomial trees, {\it Journal of
Finance} 49, no. 2, (1994), 711-818.

\medskip
\par\noindent
[S.K] Schoenmakers, J. and Kloeden, P.,
Robust option replication for a Black and Scholes model extended with
nondeterministic trends,
{\it Journal of Applied Mathematics and Stochastic Analysis} 12,
no. 2, (1999), 113-120.

\medskip
\par\noindent
[Sc] Scott, L., Option pricing when the variance changes randomly:
theory, estimation and an application, {\it Journal of Financial
and Quantitative Analysis} 22, no. 4, (1987), 419-438.

\medskip
\par\noindent
[St] Stoica, G.,
A stochastic delay financial model,
To appear in: {\it Proceedings of the American Mathematical Society} (2004).

%\newpage

\medskip
\bigskip
\noindent
Mercedes Arriojas\\
Department of Mathematics \\
University  Central of Venezuela \\
Caracas \\
Venezuela \\
e-mail: marrioja@euler.ciens.ucv.ve

\bigskip
\noindent
Yaozhong Hu \\ 
Department of Mathematics \\   
University of Kansas \\  
405 Snow Hall \\   
Lawrence \\
KS 66045-2142 \\ 
e-mail: hu@math.ukans.edu

\bigskip
\noindent
Salah-Eldin A. Mohammed \\
Department of Mathematics\\
Southern Illinois University at Carbondale\\
 Carbondale\\
 Illinois 62901, U.S.A.\\
e-mail: salah@sfde.math.siu.edu\\
web-site: http://sfde.math.siu.edu

\bigskip
\noindent
Gyula Pap \\
Institute of Mathematics and Informatics\\
University of Debrecen\\
Pf.~12, H-4010 Debrecen, Hungary\\
e-mail: papgy@inf.unideb.hu\\
web-site: http://www.inf.unideb.hu/valseg/dolgozok/papgy/papgy.html

\end{document}